\documentclass[12pt,draft,oneside]{article}
\usepackage[cp1251]{inputenc}
\usepackage[russian,ukrainian,english]{babel}
\usepackage{amsmath}
\usepackage{amsthm}
\usepackage{amssymb}
\usepackage{amsfonts}
\usepackage{latexsym}
\usepackage{calc}
\usepackage{indentfirst}
\usepackage{color}
\usepackage{cite}
\usepackage{stmaryrd}

\newcommand{\Aut}{\mathrm{Aut\mkern 2mu}}
\newcommand{\Halo}{\mathrm{Halo\mkern 2mu}}

\newtheorem{theorem}{Theorem}
\newtheorem{proposition}{Proposition}
\newtheorem{remark}{Remark}
\newtheorem{lemma}{Lemma}
\newtheorem{corollary}{Corollary}

\newcommand{\sep}{/\kern-2pt/ }

\allowdisplaybreaks

\begin{document}
	
\begin{center}
	\bf \MakeUppercase{On some classes of noncommutative dimonoids}
\end{center}
	
\begin{center}
	\bf \MakeUppercase{Gavrylkiv V.M.}
\end{center}

Department of Algebra and Geometry,  Vasyl Stefanyk Carpathian National University, Ivano-Frankivsk, Ukraine

e-mail: vgavrylkiv@gmail.com

url: https://gavrylkiv.cnu.edu.ua/

orcid: https://orcid.org/0000-0002-6256-3672

\bigskip UDC 512.53
	
\vspace{10pt plus 0.5pt} {\abstract{ \noindent V.\ M.\ Gavrylkiv.\ % Author's name (in English).
\textit{On some classes of noncommutative dimonoids} \vspace{3pt} %\English % A title of the paper (in English).
		
The present paper is devoted to the study of dimonoids, algebraic structures with  two associative binary operations that satisfy a prescribed system of axioms. We investigate the properties of dual dimonoids. In the class of noncommutative dimonoids, we construct a number of abelian, nonabelian, and rectangular dimonoids. The  structure of these objects is analyzed,  their automorphism groups and halos are computed.


}}
	
\selectlanguage{english}
	
%\allowdisplaybreaks
	
%\selectlanguage{ukrainian}
	
\vskip10pt
	
\section*{Introduction}

The notions of a dialgebra and a dimonoid were introduced by J.-L. Loday~\cite{Lod}. A {\em dimonoid} is an algebraic structure $(D,\dashv,\vdash)$ consisting of a  set $D$ equipped with two associative binary operations $\dashv$ and $\vdash$ satisfying the following axioms:
\begin{align*}
(x \dashv y) \dashv z = x \dashv (y \vdash z), \hspace{10mm}(D_1) \\
(x \vdash y) \dashv z = x \vdash (y \dashv z), \hspace{10mm} (D_2)\\
(x \dashv y) \vdash z = x \vdash (y \vdash z). \hspace{10mm} (D_3)
\end{align*}

In fact, each semigroup $(D,\dashv)$ can be consider as a dimonoid $(D,\dashv,\dashv)$, and this structure will be called a {\em trivial} dimonoid. Consequently, dimonoids are a generalization of semigroups.
Dialgebras are linear analogues of dimonoids. A {\em dialgebra} is a vector space over a field with two bilinear associative binary operations satisfying the axioms of a dimonoid and thus the results obtained for dimonoids can be applied to dialgebra theory~\cite{B,F, Lod, M}. At the present time dimonoids have became a standard tool and play a prominent role in problems from the theory of Leibniz algebras. T. Pirashvili~\cite{P} introduced the notion of a duplex which generalizes the notion of a dimonoid and constructed a free duplex. In~\cite{MDS}, a {\em $\mathfrak{g}$-dimonoid} is introduced, and the construction of the free $\mathfrak{g}$-dimonoid is described. The properties of free dimonoids were utilized in~\cite{Lod} to characterize free dialgebras and to study their cohomologies. The notion of a dimonoid was employed in~\cite{Liu} to introduce the concept of a one-sided diring and explore basic properties of dirings. Dimonoid are also associated with  restrictive  bisemigroups~\cite{Sh} and doppelsemigroups~\cite{GR1, GDS2, GUDS, GSDS, GLDS}.

The first result about dimonoids is a Loday’s description~\cite{Lod} of an absolute free dimonoid generated by a given set. Many classes of dimonoids were studied by Anatolii~Zhuchok and Yurii~Zhuchok. The independence of the axioms defining a dimonoid was established in~\cite{ZhAL2011}. A comprehensive investigation of commutative, free commutative, and free abelian dimonoids was carried out in~\cite{ZhADM2009,ZhADM2010,ZhADM2015}, respectively. Further structural aspects of dimonoids, such as dibands of subdimonoids and semilattice decompositions, were examined in~\cite{ZhMS2011,ZhDM2011}. In addition, the construction of free rectangular dimonoids, as well as free normal and free $(lr, rr)$-dibands, was presented in~\cite{ZhADM2011a,ZhADM2011,ZhADM2013}. Free abelian dibands and some of their properties were studied in~\cite{ZhVLU2017, ZhADM2018}. The least semilattice congruence on free dimonoids was explicitly described in~\cite{ZhUMJ2011}, providing important insight into the lattice-theoretic properties of these structures. The works~\cite{ZhQRL2013, ZhMP2014, ZhCA2017a, ZhIJAC2021} focused on the study of free products of dimonoids and relatively free dimonoids. The description of free left $n$-nilpotent and free left $n$-dinilpotent dimonoids was given in~\cite{ZhADM2013nil, ZhSF2016}, where special attention was paid to the behavior of nilpotency conditions within the dimonoid setting. Moreover, representations of ordered dimonoids by means of binary relations were analyzed in~\cite{ZhAEJM2014}. Significant contributions to the theory of endomorphisms and automorphisms in the context of dimonoids were made by Y. Zhuchok in~\cite{ZhBAS2014,ZhCA2017,ZhJA2024}.

In this paper, we investigate the properties of dual dimonoids, and in the class of noncommutative dimonoids we construct a number of abelian, nonabelian, and rectangular dimonoids. The internal structure of these objects is analyzed,  their automorphism groups and halos are computed.  In~\cite{Gdim2}, we use the obtained results and the results of the papers~\cite{GR1, GDS2} for classifications up to isomorphism of dimonoids with at most three elements.

\section{Preliminaries}

An element $z$ of a semigroup $(S,*)$ is called a {\em  left zero} (resp. a {\em  right zero}) in $S$ if $z*a=z$ (resp. $a*z=z$) for any $a\in S$. An element $0\in S$  is called a {\em   zero} if $0$ is a left zero and a right zero.

Let $(S,*)$ be a semigroup and $0\notin S$. The binary operation $*$ defined on  $S$  can be extended to $S\cup\{0\}$ putting $0*s=s*0=0$ for all $s\in S\cup \{0\}$.
The notation $(S,*)^{+0}$  denotes a semigroup $(S\cup\{0\},*)$ obtained from $(S,*)$ by adjoining the extra zero $0$ (regardless of whether $(S,*)$ has or has not the zero).

The semigroup is a {\em band}, if all its elements are idempotents. Commutative bands are called {\em semilattices}.

If $(S,*)$ is a semigroup, then the semigroup $(S,{*}^d)$ with operation $x{*}^d y=y* x$ is called {\em dual} to $(S,*)$, denoted $(S,*)^d$. It follows that $(S,*)^d = (S,*)$ if and only if $(S,*)$ is a commutative semigroup.

A semigroup $(S,*)$ is called a {\em null semigroup} if there exists an element $0\in S$ such that $x*y=0$ for all $x,y\in S$. In this case  $0$ is a zero of $S$. All null semigroups on the same set are isomorphic. By $O_S$ we denote a null semigroup on a set $S$. If $S$ is finite of cardinality $|S|=n$, then instead of $O_S$ we use $O_n$.  

\smallskip

Let $S$ be a nonempty set, $0\in S$ and $A\subseteq S\setminus\{0\}$. Define  the binary operation $*$ on  $S$ in the following way:

$$x* y=\begin{cases}
x, \text{ if }y=x\in A \\
0,\ \text{otherwise}.
\end{cases}$$

It is easy to check that a set $S$ endowed with the operation $*$ is a commutative semigroup with  zero $0$, and we denote this semigroup by $O^A_S$. If $A=S\setminus\{0\}$, then $O^A_S$ is a semilattice. In the case when $A$ is an emptyset,  $O^A_S$ coincides with a null  semigroup  with  zero $0$. The semigroups $O^A_S$ and $O^B_T$ are isomorphic if and only if $|S|=|T|$ and $|A|=|B|$. If $S$ is a finite set of cardinality $|S|=n$ and $|A|=m$, then we use $O_n^m$ instead of $O^A_S$.

\smallskip

A semigroup $(S,*)$ is said to be a {\em left} (resp. {\em right}) {\em zero semigroup} if $a*b=a$ (resp. $a*b=b$) for any $a,b\in S$. By  $LO_S$ and $RO_S$ we denote  a left zero semigroup and a right zero semigroup on a set $S$, respectively.  It is easy to see that the semigroups $LO_S$ and $RO_S$ are dual. If $S$ is finite of cardinality $|S|=n$, then instead of $LO_S$ and $RO_S$ we use $LO_n$ and $RO_n$, respectively.

\smallskip

Let $S$ be a nonempty set, $A\subseteq S$ and $0\notin S$. Define  the binary operation $*$ on  $S^0=S\cup\{0\}$ in the following way:

$$x* y=\begin{cases}
x, \text{ if }y \in A \\
0, \text{ if }y \in S^0\setminus A.
\end{cases}$$

It is easy to check that a set $S^0$ endowed with the operation $*$ is a semigroup with  zero $0$, and we denote this semigroup by $LO^{\sim 0}_{A\leftarrow S}$. If $A=S$, then $LO^{\sim 0}_{A\leftarrow S}$ coincides with $LO^{+0}_S$. In the case when $A$ is an emptyset,  $LO^{\sim 0}_{A\leftarrow S}$ coincides with a null semigroup $O_{S^0}$ with  zero $0$. The semigroups $LO^{\sim 0}_{A\leftarrow S}$ and $LO^{\sim 0}_{B\leftarrow T}$ are isomorphic if and only if $|S|=|T|$ and $|A|=|B|$. If $S$ is a finite set of cardinality $|S|=n$ and $|A|=m$, then we use $LO^{\sim 0}_{m\leftarrow n}$ instead of $LO^{\sim 0}_{A\leftarrow S}$.

By $RO^{\sim 0}_{A\leftarrow S}$ we denote a dual semigroup of $LO^{\sim 0}_{A\leftarrow S}$.
 
\smallskip 

Let  $a$ and $c$ be different elements of a set $S$. Define the associative binary operation $\dashv_c^a$ on  $S$ in the following way:

$$x\dashv_c^a y=\begin{cases}
a,\text{ if } x=y=a \\
c,\text{ if } x=a\text{ and } y\neq a\\
x,\text{ if } x\neq a.
\end{cases}$$

If $|S|\geq 3$, then $(S,\dashv_c^a)$ is a noncommutative band in which all elements $z\neq a$ are left zeros.

It is not difficult to check that for any different $b, d\in S$, the semigroups $(S,\dashv_c^a)$ and $(S,\dashv_d^b)$ are
isomorphic. We denote this semigroup by $LOB_S$. If $S$ is a finite set of cardinality $|S|=n$,  then we use $LOB_n$ instead of $LOB_S$. 

By $ROB_S$ we denote a dual semigroup of $LOB_S$.

\smallskip

Let $S$ be a nonempty set, $A$ be a nonempty subset of $S$, and $a\in A$. Define the associative binary operation $*$ on  $S$ in the following way:

$$x* y=\begin{cases}
x,\text{ if } x\in A \\
a,\text{ if } x \notin A.
\end{cases}$$

We denote the semigroup $(S,*)$  by $LO_{A\leftarrow S}$. It follows that  all elements $z\in A$ are left zeros of $LO_{A\leftarrow S}$. If $A=\{a\}$, then $LO_{A\leftarrow S}$ coincides with a null semigroup $O_S$ with zero $a$. If $A=S$, then $LO_{A\leftarrow S}$ coincides with a left zero semigroup $LO_S$. The semigroups $LO_{A\leftarrow S}$ and $LO_{B\leftarrow T}$ are isomorphic if and only if $|S|=|T|$ and $|A|=|B|$. If $S$ is a finite set of cardinality $|S|=n$ and $|A|=m$, then we use $LO_{m\leftarrow n}$ instead of $LO_{A\leftarrow S}$.

By $RO_{A\leftarrow S}$ we denote a dual semigroup of $LO_{A\leftarrow S}$.

\bigskip

An element $e$ of a dimonoid $(D,\dashv, \vdash)$ is called a {\em bar-unit}~\cite{Lod} if $e \vdash d = d = d \dashv e$ for all $d\in D$. In contrast to monoids a dimonoid may have many bar-units. The set of all bar-units of a dimonoid $(D,\dashv, \vdash)$ is called the {\em halo} of $(D,\dashv, \vdash)$~\cite{Lod}, denoted $\Halo(D,\dashv, \vdash)$. A nonempty subset $B\subseteq D$ is called a {\em subdimonoid} of a dimonoid $(D,\dashv, \vdash)$  if $a \dashv b, a \vdash b \in B$ for any $a, b \in B$. If the halo of $(D,\dashv, \vdash)$ is nonempty, then it is a subdimonoid of $(D,\dashv, \vdash)$. 

An element $0\in D$  is called  a {\em   zero} of a dimonoid $(D,\dashv, \vdash)$ if $0$ is a zero of $(D,\dashv)$ and a  zero of $(D,\vdash)$. Let $(D,\dashv, \vdash)$ be a dimonoid and $0\notin D$. The binary operations defined on  $D$  can be extended to $D\cup\{0\}$ putting $0\dashv d=d\dashv 0=0=0\vdash d=d\vdash 0 $ for all $d\in D\cup \{0\}$. The notation $(D,\dashv, \vdash)^{+0}$  denotes a dimonoid $D\cup\{0\}$ obtained from $D$ by adjoining the extra zero $0$. It follows  that $\Halo((D,\dashv, \vdash)^{+0}) = \Halo(D,\dashv, \vdash)$.

A bijective map $\psi : D_1 \to D_2$ is called an {\em isomorphism } from a dimonoid $(D_1,\dashv_1, \vdash_1)$ to a dimonoid $(D_2,\dashv_2, \vdash_2)$ if $$\psi(a\dashv_1 b)=\psi(a)\dashv_2\psi(b)\ \ \text{  and  }\ \ \psi(a\vdash_1 b)=\psi(a)\vdash_2\psi(b)$$ for all $a,b\in D_1$.

If there exists an isomorphism from a dimonoid $(D_1,\dashv_1, \vdash_1)$ to a dimonoid $(D_2,\dashv_2, \vdash_2)$, then $(D_1, \dashv_1, \vdash_1)$ and $(D_2, \dashv_2, \vdash_2)$ are said to be {\em isomorphic}, denoted $(D_1,\dashv_1, \vdash_1)\cong (D_2,\dashv_2, \vdash_2)$. An isomorphism $\psi: D\to D$ is called an {\em   automorphism} of a dimonoid $(D,\dashv, \vdash)$. By $\Aut(D,\dashv, \vdash)$ we denote the automorphism group of a dimonoid $(D,\dashv, \vdash)$. It follows that $\Aut((D,\dashv, \vdash)^{+0}) = \Aut(D,\dashv, \vdash)$.

For a dimonoid $(D,\dashv,\vdash)$, if $\mathbb S$ and $\mathbb T$ denote the semigroups $(D,\dashv)$ and $(D,\vdash)$, respectively, then $\mathbb S \rbag \mathbb T$ stands for the dimonoid $(D,\dashv,\vdash)$.

\smallskip

Following the algebraic tradition, for a set $X$ by $S_X$ we denote the group of all bijections of $X$.

\section{Properties of dual dimonoids}\label{sec:dual}

This section considers  properties of dual dimonoids and certain other basic properties that will be utilized in the forthcoming sections.

Let $(D,\dashv, \vdash)$ be a dimonoid. Define new operations $\dashv^d$ and  $\vdash^d$ on $D$ by 
$$x \dashv^d y = y \vdash x\ \ \text{  and  }\ \  x \vdash^d y = y \dashv x.$$
It is easy to verify that $(D,\dashv^d, \vdash^d)$ forms a new dimonoid, called the {\em dual dimonoid} of $(D,\dashv, \vdash)$~\cite{Lod}, and denoted by $(D,\dashv, \vdash)^d$. 
It follows that the unary duality operation is involutive in the sense that $((D,\dashv, \vdash)^d)^d=(D,\dashv, \vdash)$. In fact, $(D,\dashv, \vdash)^d$ is a dimonoid if and only if $(D,\dashv, \vdash)$ is a dimonoid.  As usual, we say that a dimonoid $(D,\dashv, \vdash)$ is {\em self-dual} if $(D,\dashv, \vdash)^d=(D,\dashv, \vdash)$.

\begin{remark}Observe that if we put
$x \dashv' y = y \dashv x$  and  $x \vdash' y = y \vdash x$, then in general case $(D,\dashv', \vdash')$ is not a dimonoid. Indeed, let us consider a left zero and a right zero dimonoid $LO_D\rbag RO_D = (D, \dashv, \vdash)$ with operations $x\dashv y = x$ and $x\vdash y=y$, see~\cite{Lod}. Taking into account that $x \dashv' y = y \dashv x = y$ and $x \vdash' y = y \vdash x = x$, we conclude that $ (x \dashv' y) \dashv' z = z$ while   $x \dashv' (y \vdash' z) = x \dashv' y = y$  for all $x,y,z\in D$. Therefore, the axiom $(D_1)$ does not hold, and  $(D,\dashv', \vdash')$ can not be a dimonoid.
\end{remark}

A dimonoid $(D,\dashv,\vdash)$ is called {\em commutative}~\cite{ZhADM2009} if both semigroups $(D,\dashv)$ and $(D,\vdash)$ are commutative.

\begin{proposition}\label{duality} Let $(D,\dashv, \vdash)$ be a dimonoid. Then $\Aut(D,\dashv, \vdash)^d = \Aut(D,\dashv, \vdash)$ and $\Halo(D,\dashv, \vdash)^d = \Halo(D,\dashv, \vdash)$. Moreover, $(D,\dashv, \vdash)$ is commutative if and only if $(D,\dashv, \vdash)^d$ is commutative.
\end{proposition}

\begin{proof} Let $\psi: D\to D$ be an automorphism of a dimonoid $(D,\dashv, \vdash)$. Then $\psi$ is a bijection. Taking into account that $\psi(x\dashv^d y) = \psi(y\vdash x) = \psi(y)\vdash \psi(x) = \psi(x) \dashv^d \psi(y)$ and
$\psi(x\vdash^d y) = \psi(y\dashv x) = \psi(y)\dashv \psi(x) = \psi(x) \vdash^d \psi(y)$ for all $x,y\in D$, we conclude that $\psi$ is an automorphism of a dimonoid $(D,\dashv, \vdash)^d$ as well. Using involutivity of the unary duality operation, we conclude that $\Aut(D,\dashv, \vdash)^d = \Aut(D,\dashv, \vdash)$.

Let $e$ be a bar-unit of a dimonoid $(D,\dashv, \vdash)$. Then $e\vdash d = d = d\dashv e$ for all $d\in D$. Since $e\vdash^d d = d\dashv e = d = e\vdash d = d\dashv^d e$ for all $d\in D$, it follows that $e$ is a bar-unit of a dimonoid $(D,\dashv, \vdash)^d$ as well. Using involutivity of the unary duality operation, we conclude that $\Halo(D,\dashv, \vdash)^d = \Halo(D,\dashv, \vdash)$.

Since commutativity of an operation $\dashv$ is equivalent to commutativity of an operation $\vdash^d$ and commutativity of an operation $\vdash$ is equivalent to commutativity of an operation $\dashv^d$, we conclude that $(D,\dashv, \vdash)$ is commutative if and only if $(D,\dashv, \vdash)^d$ is commutative.
\end{proof}

A dimonoid $(D,\dashv, \vdash)$ is called {\em abelian}~\cite{ZhADM2015} if $x \dashv y = y \vdash x$ for all $x,y\in D$.  

\begin{proposition}\label{da}
Let $(D,\dashv, \vdash)$ be a dimonoid. Then the following conditions are equivalent:
\begin{itemize}
\item[1)] $(D,\dashv)$ and $(D,\vdash)$ are dual semigroups;
\item[2)] $(D,\dashv, \vdash)$ is abelian;
\item[3)] $(D,\dashv, \vdash)$ is self-dual.
\end{itemize}
\end{proposition}

\begin{proof} $(1)\Rightarrow(2)$ If $(D,\dashv)$ and $(D,\vdash)$ are dual semigroups, then $x\dashv y = y\vdash x $ for all $x,y\in D$, and hence $(D,\dashv, \vdash)$ is an abelian dimonoid.

$(2)\Rightarrow(3)$ Let $(D,\dashv, \vdash)$ be an abelian dimonoid. Taking into account that 
$$x \dashv^d y = y \vdash x = x \dashv y\ \ \text{  and  }\ \  x \vdash^d y = y \dashv x = x \vdash y,$$
we conclude that $\dashv^d\ =\ \dashv$ and $\vdash^d\ =\ \vdash$, and hence $(D,\dashv, \vdash)$ is a self-dual dimonoid.

$(3)\Rightarrow(1)$ If $(D,\dashv, \vdash)$ is a self-dual dimonoid, then $x \dashv y = x \dashv^d y = y \vdash x$ for all $x,y\in D$, and hence $(D,\dashv)$ and $(D,\vdash)$ are dual semigroups.
\end{proof}

\begin{corollary} Nonabelian dimonoids are divided into the pairs of dual dimonoids.
\end{corollary}

Taking into account that automorphism groups of dual semigroups coincide, we conclude the following corollary.

\begin{corollary}\label{autabdim} Let $(D,\dashv, \vdash)$ be an abelian dimonoid. Then $\Aut(D,\dashv, \vdash) = \Aut(D,\dashv) = \Aut(D,\vdash)$.
\end{corollary}

The definition of the halo of a dimonoid $(D,\dashv, \vdash)$ and Proposition~\ref{da} implies the following corollary.

\begin{corollary}\label{haloabdim} Let $(D,\dashv, \vdash)$ be an abelian dimonoid. If the halo of $(D,\dashv, \vdash)$ is nonempty, then it coincides with the subsemigroup of all right identities of a semigroup $(D,\dashv)$ and the subsemigroup of all left identities of a semigroup $(D,\vdash)$.
\end{corollary}

Since commutative semigroups $(D,\dashv)$ and $(D,\vdash)$ are dual if and only if their operations coincide, Proposition~\ref{da} implies the following corollary.

\begin{corollary} Commutative nontrivial dimonoids are nonabelian.
\end{corollary}

On the other hand, it is clear to see that all commutative trivial dimonoids are abelian and all noncommutative trivial dimonoids are nonabelian. A  left zero and a right zero dimonoid $LO_D\rbag RO_D$ is an example of a nontrivial abelian noncommuative dimonoid. In~\cite{Gdim2}, we give examples of commuative nonabelian dimonoids, see also~\cite{ZhADM2009}.

A semigroup $(S, *)$ is called {\em right}  {\em commutative}~\cite{ZhADM2013}, if it satisfies the identity $s*x*y = s*y*x$  for all $s, x, y \in S$. It follows that if a semigroup $(S, *)$ is right commutative, then a semigroup $(S, *)^{+0}$ is right commutative as well.

The following lemma was proved by A.~Zhuchok~\cite{ZhADM2013}.

\begin{lemma}[A. Zhuchok]\label{rcommut}
Let $(D, \dashv)$ be an arbitrary semigroup and $(D, \vdash)$ be a dual semigroup to $(D, \dashv)$. Then an algebraic structure $(D, \dashv, \vdash)$ is a dimonoid if and only if $(D, \dashv)$ is a right commutative semigroup.
\end{lemma}

It follows that the task of describing abelian dimonoids reduces to the task of describing right commutative semigroups.

\section{Abelian noncommutative dimonoids}\label{sec:abelian}

In this section we construct new examples of abelian noncommutative dimonoids using some right commutative semigroups and also calculate their halos and automorphism groups.

\begin{proposition}\label{rcommLOAD}  Let $A$ be a nonempty proper subset of a set $D$, and $a\in A$. Then  $LO_{A\leftarrow D}=(D,\dashv)$, where
$$x\dashv y=\begin{cases}
x,\text{ if } x\in A \\
a,\text{ if } x \in D\setminus A,
\end{cases}$$
is a right commutative semigroup with the automorphism group $\Aut(LO_{A\leftarrow D})=S_{A\setminus\{a\}}\times S_{D\setminus A}$.
\end{proposition}
\begin{proof}It is immediate to check that  $LO_{A\leftarrow D}$ is a right commutative semigroup. 
%If $A=D$ the semigroup $LO_{A\leftarrow D}= LO_D$ coincide with a left zero semigroup on $D$, and hence in this case $\Aut(LO_{A\leftarrow D})=S_D$.  

Let $\psi$ be an arbitrary automorphism of $LO_{A\leftarrow D}$. Since automorphisms preserves left zeros, it follows that $\psi(A)=A$, and thus $\psi(D\setminus A)=D\setminus A$. Fix any $d\in D\setminus A$. It follows that $\psi(d)\in D\setminus A$, and hence $\psi(a)=\psi(d\dashv a) = \psi(d)\dashv \psi(a)=a$.

On the other hand, let $f$ be any bijection of $D$ such that $f(A)=A$, $f(D\setminus A)=D\setminus A$ and $f(a)=a$.  If $x\in A$, then $f(x)\in A$, and hence $f(x\dashv y)= f(x)= f(x)\dashv f(y)$ for all $y\in D$.
In the case $x\in D\setminus A$ we have that $f(x)\in D\setminus A$, and thus $f(x\dashv y)= f(a)= a = f(x)\dashv f(y)$ for all $y\in D$.
It follows that any bijection of $A$ that preserves $a$ and any bijection of $D\setminus A$ generate an automorphism of $LO_{A\leftarrow D}$. Therefore, $\Aut(LO_{A\leftarrow D})=S_{A\setminus\{a\}}\times S_{D\setminus A}$.
\end{proof}

Lemma~\ref{rcommut}, Propositions~\ref{rcommLOAD} and~\ref{da}, Corollaries~\ref{autabdim} and~\ref{haloabdim} and a definition of $LO_{A\leftarrow D}$ imply the following theorem.

\begin{theorem}\label{rcommLOADdim} Let $A$ be a nonempty proper subset of a set $D$, and $a\in A$. Then  $LO_{A\leftarrow D}\rbag RO_{A\leftarrow D}$ is a abelian dimonoid with empty halo and $\Aut(LO_{A\leftarrow D}\rbag RO_{A\leftarrow D})=S_{A\setminus\{a\}}\times S_{D\setminus A}$. Moreover, for $|A|>1$ the dimonoid $LO_{A\leftarrow D}\rbag RO_{A\leftarrow D}$ is not commutative.
\end{theorem}

\begin{proposition}\label{rcommLOAD0} Let $D$ be a nonempty set, $A\subseteq D$, $0\notin D$ and  $D^0=D\cup\{0\}$. Then  $LO^{\sim 0}_{A\leftarrow D}=(D^0,\dashv)$, where
$$x\dashv y=\begin{cases}
x, \text{ if }y \in A \\
0, \text{ if }y \in D^0\setminus A.
\end{cases}$$
is a right commutative semigroup with the automorphism group $\Aut(LO^{\sim 0}_{A\leftarrow D})=S_A\times S_{D\setminus A}$.
\end{proposition}
\begin{proof}

Since all elements of a set $A$ are right identities of $LO^{\sim 0}_{A\leftarrow D}$, we conclude that $s\dashv x\dashv y = s = s\dashv y\dashv x$ for all $x, y \in A$, $s\in D^0$. If $x\notin A$ or $y\notin A$, then $s\dashv x\dashv y = 0 = s\dashv y\dashv x$ for all $s, x, y \in D^0$.  Therefore, $LO^{\sim 0}_{A\leftarrow D}$ is a right commutative semigroup.

Let $\psi$ be an arbitrary automorphism of $LO^{\sim 0}_{A\leftarrow D}$. It follows that $0$ is a zero and $A$ is a subsemigroup of all right identities of the semigroup $LO^{\sim 0}_{A\leftarrow D}$. Since automorphisms preserves a zero and right identities, it follows that $\psi(0)= 0$, $\psi(A)=A$, and thus $\psi(D\setminus A)=D\setminus A$. 

On the other hand, let $f$ be any bijection of $D^0$ such that $f(0)=0$, $f(A)=A$, and hence $f(D\setminus A)=D\setminus A$.  If $y\in A$, then $f(y)\in A$, and hence $f(x\dashv y)= f(x)= f(x)\dashv f(y)$ for all $x\in D^0$.
In the case $y\in D^0\setminus A$ we have that $f(y)\in D^0\setminus A$, and thus $f(x\dashv y)= f(0)= 0 = f(x)\dashv f(y)$ for all $x\in D^0$.
It follows that any bijection of $A$  and any bijection of $D\setminus A$ generate an automorphism of $LO^{\sim 0}_{A\leftarrow D}$. Therefore, $\Aut(LO^{\sim 0}_{A\leftarrow D})=S_A\times S_{D\setminus A}$.
\end{proof}

Lemma~\ref{rcommut}, Propositions~\ref{rcommLOAD0} and~\ref{da}, Corollaries~\ref{autabdim} and~\ref{haloabdim} and a definition of $LO^{\sim 0}_{A\leftarrow D}$ imply the following theorem.

\begin{theorem}\label{rcommLOAD0dim} Let $D$ be a nonempty set and $A\subseteq D$. Then  $LO^{\sim 0}_{A\leftarrow D}\rbag RO^{\sim 0}_{A\leftarrow D}$ is and abelian dimonoid with $\Halo(LO^{\sim 0}_{A\leftarrow D}\rbag RO^{\sim 0}_{A\leftarrow D}) = A$ and $\Aut(LO^{\sim 0}_{A\leftarrow D}\rbag RO^{\sim 0}_{A\leftarrow D})=S_A\times S_{D\setminus A}$. Moreover, if $A$ is a nonempty set and $|D|>1$, then the dimonoid $LO^{\sim 0}_{A\leftarrow D}\rbag RO^{\sim 0}_{A\leftarrow D}$ is not commutative.
\end{theorem}

\begin{proposition}\label{rcommLOB} Let $D$ be a set of cardinality $|D|\geq 2$ and $a$, $c$ be different elements of $D$. Then  $LOB_D=(D,\dashv)$, where
$$x\dashv y=\begin{cases}
a,\text{ if } x=y=a \\
c,\text{ if } x=a\text{ and } y\neq a\\
x,\text{ if } x\neq a.
\end{cases}$$
is a right commutative semigroup with the automorphism group $\Aut(LOB_D)=S_{D\setminus\{a,c\}}$.
\end{proposition}
\begin{proof}
If $s$ is a left zero, then immediate $s\dashv x\dashv y = s = s\dashv y\dashv x$ for all $x, y \in D$. Taking into account that $a\dashv x \dashv y\in\{a,c\}$ for all $x, y \in D$ and in the semigroup $LOB_D$ the equation $a\dashv u \dashv v = a$ has a unique solution $u=v=a$, we conclude that $a\dashv x\dashv y = a\dashv y\dashv x$ for all $x, y \in D$. Therefore, $LOB_D$ is a right commutative semigroup.

Let $\psi$ be an arbitrary automorphism of $LOB_D$. Since $a$ is a unique right identity of $LOB_D$ and automorphisms preserves right identities, it follows that $\psi(a)= a$. Fix any $d\ne a$. It follows that $\psi(d)\ne a$, and hence $\psi(c)=\psi(a\dashv d) = \psi(a)\dashv \psi(d)= a\dashv \psi(d) = c$.

On the other hand, let $f$ be any bijection of $D$ such that $f(a)= a$ and $f(c)=c$.  If $x\ne a$, then $f(x)\ne a$, and hence $f(x\dashv y)= f(x)= f(x)\dashv f(y)$ for all $x, y\in D$. It follows that $f(a\dashv a)= f(a)= a = a\dashv a = f(a)\dashv f(a)$.
In the case $y\ne a$ we have that $f(y)\ne a$, and thus $f(a\dashv y)= f(c)= c = a\dashv f(y) = f(a)\dashv f(y)$ for all $y\in D$.
Therefore, any bijection of $D$ that preserves $a$ and $c$ generates an automorphism of $LOB_D$. Consequently, $\Aut(LOB_D)=S_{D\setminus\{a,c\}}$.
\end{proof}

Lemma~\ref{rcommut}, Propositions~\ref{rcommLOB} and~\ref{da}, Corollaries~\ref{autabdim} and~\ref{haloabdim} and a definition of $LOB_D$ imply the following theorem.

\begin{theorem}\label{rcommLOBdim} Let $D$ be a set of cardinality $|D|\geq 2$ and  $a$, $c$ be different elements of $D$. Then $LOB_D\rbag ROB_D$ is an abelian dimonoid with $\Halo(LOB_D\rbag ROB_D) = \{a\}$ and $\Aut(LOB_D\rbag ROB_D)=S_{D\setminus\{a,c\}}$. Moreover, if $|D|\geq 3$, then $LOB_D\rbag ROB_D$ is a noncommutative dimonoid.
\end{theorem}

\begin{proposition}\label{rcommLO+0} Let $D$ be a nonempty set. Then  $LO_D^{+0}$ is a right commutative semigroup with the automorphism group $\Aut(LO_D^{+0})=S_D$.
\end{proposition}
\begin{proof} Since $LO_D$ is a right commutative semigroup, $LO_D^{+0}$ is a right commutative semigroup as well. It follows that $\Aut(LO_D^{+0}) = \Aut(LO_D)= S_D$.
\end{proof}

Lemma~\ref{rcommut}, Propositions~\ref{rcommLO+0} and~\ref{da}, Corollaries~\ref{autabdim} and~\ref{haloabdim} and a definition of $LO_D^{+0}$ imply the following theorem.

\begin{theorem}\label{rcommLO+0dim} Let $D$ be a nonempty set. Then  $(LO_D\rbag RO_D)^{+0} = LO_D^{+0}\rbag RO_D^{+0}$ is an abelian dimonoid with $\Halo((LO_D\rbag RO_D)^{+0}) = D$ and $\Aut((LO_D\rbag RO_D)^{+0})=S_D$. Moreover, if $|D|>1$, then $(LO_D\rbag RO_D)^{+0}$ is a noncommutative  dimonoid.
\end{theorem}

\section{Nonabelian noncommutative dimonoids}\label{sec:nonabelian}

In this section we construct some new examples of nonabelian noncommutative dimonoids with empty halos and calculate their automorphism groups. 

\begin{theorem}
Let $D$ be a set of cardinality $|D|>2$, and $a, c$ be different elements of $D$. An algebraic structure $LOB_D \rbag O_D^{\{a\}}=(D, \dashv, \vdash)$, where

\begin{center}
$x\dashv y=\begin{cases}
a,\text{ if } x=y=a \\
c,\text{ if } x=a\text{ and } y\neq a\\
x,\text{ if } x\neq a
\end{cases}$
and\ \ \ \  
$x\vdash y=\begin{cases}
a,\text{ if } x=y=a \\
c,\ \text{otherwise}
\end{cases}$
\end{center}
is a nonabelian noncommutative dimonoid with empty halo and $\Aut(LOB_D \rbag O_D^{\{a\}})=S_{D\setminus\{a,c\}}$.

\end{theorem}

\begin{proof} Let us show that $(D, \dashv, \vdash)$ is a dimonoid dividing into cases of axioms of a dimonoid. 

{\noindent \bf Case $(D_1)$}. If $x\ne a$, then $x$ is a left zero of a semigroup $(D, \dashv)$, and hence
$(x \dashv y) \dashv z = x \dashv (y \dashv z)= x = x \dashv (y \vdash z)$ for any $y,z\in D$.
If $x=a$, then $(x \dashv y) \dashv z = x \dashv (y \dashv z), x \dashv (y \vdash z)\in\{a,c\}$ for any $y,z\in D$. Since $(x \dashv y) \dashv z=a$ if and only if $x=y=z=a$, $x \dashv (y \vdash z)=a$ if and only if $x=y=z=a$, we conclude that $(D_1)$ holds.

{\noindent \bf Case $(D_2)$}. It follows that $x \vdash (y \dashv z)\in \{a,c\}$ for any $x,y,z\in D$. Taking into account that $x \vdash y\in \{a,c\}$ for any $x,y\in D$, we conclude that $(x \vdash y) \dashv z\in \{a,c\}$ for any $x,y,z\in D$ as well. Since $(x \vdash y) \dashv z=a$ if and only if $x=y=z=a$, $x \vdash (y \dashv z)=a$ if and only if $x=y=z=a$, we conclude that $(D_2)$ holds. 

{\noindent \bf Case $(D_3)$}. Taking into account that $(x \dashv y) \vdash z, x \vdash (y \vdash z)\in \{a,c\}$ for any $x,y,z\in D$ and  $(x \dashv y) \vdash z=a$ if and only if $x=y=z=a$, $x \vdash (y \vdash z)=a$ if and only if $x=y=z=a$, we conclude that $(D_3)$ holds.

\smallskip

 Let $z\in D\setminus\{a,c\}$. Taking into account that $z\dashv a =z$ but $a\vdash z = c\ne z$, we conclude that $(D, \dashv, \vdash)$ is a nonabelian dimonoid. Since $(D, \dashv)$ is a noncommutative semigroup, $(D, \dashv, \vdash)$ is a noncommutative dimonoid.
 
Since the commutative semigroup $O_D^{\{a\}}$ contains no identities, the halo of the dimonoid $LOB_D \rbag O_D^{\{a\}}$ is emptyset.

Let $\psi$ be an arbitrary automorphism of the dimonoid $LOB_D \rbag O_D^{\{a\}}$. Then $\psi$ is an automorphism of the semigroup $LOB_D$, and hence $\psi(a)= a$, $\psi(c)= c$ according to proof of Proposition~\ref{rcommLOB}.

On the other hand, let $f$ be any bijection of $D$ such that $f(a)= a$ and $f(c)=c$.  Then $f$ is an automorphism of the semigroup $LOB_D$ according to proof of Proposition~\ref{rcommLOB}.

If $x\ne a$ or $y\ne a$, then $f(x)\ne a$ or $f(y)\ne a$, and hence $f(x\vdash y)= f(c)= c = f(x)\vdash f(y)$ for all $x, y\in D$. Taking into account that $f(a\vdash a)= f(a)= a = a\vdash a = f(a)\vdash f(a)$, we conclude that
$f$ is an automorphism of the semigroup $O_D^{\{a\}}$, and hence $f$ is an automorphism of the dimonoid $LOB_D \rbag O_D^{\{a\}}$.

Therefore, any bijection of $D$ that preserves $a$ and $c$ generates an automorphism of the dimonoid $LOB_D \rbag O_D^{\{a\}}$. Consequently, $\Aut(LOB_D \rbag O_D^{\{a\}})=S_{D\setminus\{a,c\}}$.
\end{proof}

\begin{theorem}
Let $D$ be a set of cardinality $|D|>1$, $a\in D$, $0\notin D$ and  $D^0=D\cup\{0\}$. An algebraic structure $LO^{\sim 0}_{{\{a\}}\leftarrow D}\rbag O^{\{a\}}_{D^{0}} =(D^0, \dashv, \vdash)$, where

\begin{center}
$x\dashv y=\begin{cases}
x,\text{ if } y=a \\
0,\text{ if } y\neq a
\end{cases}$
and\ \ \ \  
$x\vdash y=\begin{cases}
a,\text{ if } x=y=a \\
0,\ \text{otherwise}
\end{cases}$
\end{center}
is a nonabelian noncommutative dimonoid with empty halo and $\Aut(LO^{\sim 0}_{{\{a\}}\leftarrow D}\rbag O^{\{a\}}_{D^{0}})=S_{D\setminus\{a\}}$.

\end{theorem}

\begin{proof} Let us show that $(D^0, \dashv, \vdash)$ is a dimonoid dividing into cases of axioms of a dimonoid. 

{\noindent \bf Case $(D_1)$}. Let $y=z=a$, then $y\dashv z = a = y\vdash z$. Since $a$ is a right identity of a semigroup $(D^0, \dashv)$, we conclude that in this case $(x \dashv y) \dashv z = x \dashv (y \dashv z)= x \dashv a = x$ and $x \dashv (y \vdash z)=x \dashv a = x$ for any $x\in D^0$.
If $y\ne a$ or $z\ne a$, then $y\dashv z \ne a$ and $y\vdash z\ne a$. It follows that $(x \dashv y) \dashv z = x \dashv (y \dashv z)= 0$ and $x \dashv (y \vdash z)= 0$ for any $x, y, z\in D^0$.

{\noindent \bf Case $(D_2)$}. It follows that $x \vdash (y \dashv z)\in \{a,0\}$ for any $x,y,z\in D^0$. Taking into account that $x \vdash y\in \{a,0\}$ for any $x,y\in D^0$, we conclude that $(x \vdash y) \dashv z\in \{a,0\}$ for any $x,y,z\in D^0$ as well. Since $(x \vdash y) \dashv z=a$ if and only if $x=y=z=a$, $x \vdash (y \dashv z)=a$ if and only if $x=y=z=a$, we conclude that $(D_2)$ holds. 

{\noindent \bf Case $(D_3)$}. Taking into account that $(x \dashv y) \vdash z, x \vdash (y \vdash z)\in \{a,0\}$ for any $x,y,z\in D$ and  $(x \dashv y) \vdash z=a$ if and only if $x=y=z=a$, $x \vdash (y \vdash z)=a$ if and only if $x=y=z=a$, we conclude that $(D_3)$ holds.

\smallskip

 Let $z\in D\setminus\{a\}$. Taking into account that $z\dashv a =z$ but $a\vdash z = 0\ne z$, we conclude that $(D^0, \dashv, \vdash)$ is a nonabelian dimonoid. Since$(D^0, \dashv)$ is a noncommutative semigroup, $(D^0, \dashv, \vdash)$ is a noncommutative dimonoid.
 
Since the commutative semigroup $O^{\{a\}}_{D^{0}}$ contains no identities, the halo of the dimonoid $LO^{\sim 0}_{{\{a\}}\leftarrow D} \rbag O^{\{a\}}_{D^{0}}$ is emptyset.

Let $\psi$ be an arbitrary automorphism of the dimonoid $LO^{\sim 0}_{{\{a\}}\leftarrow D} \rbag O^{\{a\}}_{D^{0}}$. Then $\psi$ is an automorphism of the semigroup $LO^{\sim 0}_{{\{a\}}\leftarrow D}$, and hence $\psi(0)= 0$, $\psi(a)= a$  according to proof of Proposition~\ref{rcommLOAD0}.

On the other hand, let $f$ be any bijection of $D$ such that $f(0)= 0$ and $f(a)=a$.  Then $f$ is an automorphism of the semigroup $LO^{\sim 0}_{{\{a\}}\leftarrow D}$ according to proof of Proposition~\ref{rcommLOAD0}.

If $x\ne a$ or $y\ne a$, then $f(x)\ne a$ or $f(y)\ne a$, and hence $f(x\vdash y)= f(0)= 0 = f(x)\vdash f(y)$ for all $x, y\in D^0$. Taking into account that $f(a\vdash a)= f(a)= a = a\vdash a = f(a)\vdash f(a)$, we conclude that
$f$ is an automorphism of the semigroup $O^{\{a\}}_{D^{0}}$, and hence $f$ is an automorphism of the dimonoid $LO^{\sim 0}_{{\{a\}}\leftarrow D} \rbag O^{\{a\}}_{D^{0}}$.

Therefore, any bijection of $D^0$ that preserves $0$ and $a$ generates an automorphism of the dimonoid $LO^{\sim 0}_{{\{a\}}\leftarrow D} \rbag O^{\{a\}}_{D^{0}}$. Consequently, $\Aut(LO^{\sim 0}_{{\{a\}}\leftarrow D} \rbag O^{\{a\}}_{D^{0}})=S_{D\setminus\{a\}}$.
\end{proof}

\section{Rectangular dimonoids}\label{sec:rectangular}

A semigroup $(S,*)$ is called {\em rectangular}~\cite{ZhCA2017a} if $x*y*z = x*z$ foll all $x, y, z \in S$. A dimonoid $(D,\dashv, \vdash)$ is said to be {\em rectangular} if both $(D,\dashv)$ and $(D,\vdash)$ are rectangular semigroups. Well-known examples of rectangular dimonoids are dimonoids $LO_D\rbag RO_D$, $O_D\rbag RO_D$ and $LO_D\rbag O_D$, see~\cite{Lod, ZhADM2009, ZhADM2011a}.

It is immediate to prove the following proposition using the definition of a dimonoid, see also~\cite{ZhCA2017a}.

\begin{proposition}\label{lrecdim}
Let  $(D,\dashv)$ be a left zero semigroup and $(D, \vdash)$ be an arbitrary semigroup. An algebraic structure $(D,\dashv, \vdash)$ is a dimonoid if and only if $(D, \vdash)$ is a rectangular semigroup.
\end{proposition}

\begin{theorem} Let $A$ be a nonempty subset of a set $D$. Then $LO_D\rbag LO_{A\leftarrow D}$ and $LO_D\rbag RO_{A\leftarrow D}$ are rectangular dimonoids, moreover if $|D|>1$, then  $LO_D\rbag LO_{A\leftarrow D}$ is a nonabelian noncommutative dimonoid. If $A$ is a proper subset of a set $D$, then  $LO_D\rbag RO_{A\leftarrow D}$ is  nonabelian noncommutative as well.
\end{theorem}

\begin{proof}
It is immediate to check that  $LO_{A\leftarrow D}$ and $RO_{A\leftarrow D}$ are rectangular semigroups. Therefore, Proposition~\ref{lrecdim} implies that $LO_D\rbag LO_{A\leftarrow D}$ and $LO_D\rbag RO_{A\leftarrow D}$ are rectangular dimonoids. Since in the case $|D|>1$ the semigroups $LO_D$ and $LO_{A\leftarrow D}$ are not dual, by Proposition~\ref{da} the dimonoid $LO_D\rbag LO_{A\leftarrow D}$ is not abelian. If $|D|>1$, then  $LO_D$ is not a commutative semigroup, and hence the dimonoid $LO_D\rbag LO_{A\leftarrow D}$ is not commutative. Taking into account that for a proper subset $A\subset D$  the semigroups $LO_D$ and $RO_{A\leftarrow D}$ are not dual and the semigroup $LO_D$ is not commutative, we conclude that $LO_D\rbag RO_{A\leftarrow D}$ is a nonabelian noncommutative dimonoid as well according to Proposition~\ref{da}.
\end{proof}

\begin{proposition} Let $A$ be a nonempty proper subset of a set $D$, and $a\in A$. Then $\Halo(LO_D\rbag LO_{A\leftarrow D}) = \Halo(LO_D\rbag RO_{A\leftarrow D}) = \emptyset$ and $\Aut(LO_D\rbag LO_{A\leftarrow D}) = \Aut(LO_D\rbag RO_{A\leftarrow D}) = S_{A\setminus\{a\}}\times S_{D\setminus A}$. 
\end{proposition}
\begin{proof} Since for a proper subset $A\subset D$ the semigroups $LO_{A\leftarrow D}$ and $RO_{A\leftarrow D}$ contain no a left identity, $\Halo(LO_D\rbag LO_{A\leftarrow D}) = \Halo(LO_D\rbag RO_{A\leftarrow D}) = \emptyset$. Taking into account that $\Aut(LO_D)=S_D$ and the automorphism groups of dual semigroups coincide, we conclude that
$\Aut(LO_D\rbag LO_{A\leftarrow D}) = \Aut(LO_{A\leftarrow D}) = S_{A\setminus\{a\}}\times S_{D\setminus A}$ and
$\Aut(LO_D\rbag RO_{A\leftarrow D}) = \Aut(RO_{A\leftarrow D}) = \Aut(LO_{A\leftarrow D}) = S_{A\setminus\{a\}}\times S_{D\setminus A}$ according to Proposition~\ref{rcommLOAD}.
\end{proof}

\begin{proposition}\label{lodim}
Let $(D, \dashv)$ be a semigroup  and $(D,\vdash)$ be a null semigroup with zero $0$. An algebraic structure $(D,\dashv, \vdash)$ is a dimonoid if and only if $0$ is a left zero of $(D, \dashv)$ and $x\dashv y\dashv z = x\dashv 0$ for all $x,y,z\in D$.
\end{proposition}
\begin{proof}
Since $(D,\vdash)$ is a null semigroup,  the axiom $(D_3)$  holds:
$$(x \dashv y) \vdash z = 0 = x \vdash (y \vdash z).$$
  
Taking into account that  $x \dashv ( y \vdash z) = x \dashv 0$ for any $x,y,z\in D$, we conclude that the axiom $(D_1)$ holds if and only if $x \dashv y \dashv z = x \dashv 0$ for any $x,y,z\in D$.  Since $x \vdash (y \dashv z) = 0$ and $(x \vdash y) \dashv z = 0\dashv z$ for any $x,y,z\in D$, the axiom $(D_2)$ holds if and only if $0\dashv z = 0$ for all $z\in D$ if and only if $0$  is a left zero of $(D, \dashv)$.
\end{proof}

\begin{theorem} Let $A$ be a  subset of a set $D$ and $0\in A$. Then $LO_{A\leftarrow D}\rbag  O_D = (D,\dashv, \vdash)$, where $(D,\vdash)$ is a null semigroup with zero $0$ and
\begin{center}
$x\dashv y=\begin{cases}
x,\text{ if } x\in A \\
0,\text{ if } x \notin A,
\end{cases}$
\end{center}
is a rectangular dimonoid. Moreover if $|A|>1$, then  $LO_{A\leftarrow D}\rbag  O_D$ is a nonabelian noncommutative dimonoid. 
\end{theorem}

\begin{proof}
Using the definition of the semigroup $LO_{A\leftarrow D}$, it is immediate to check that $x\dashv y\dashv z = x\dashv 0$ for all $x,y,z\in D$. Since $0$ is a left zero of $LO_{A\leftarrow D}$, Proposition~\ref{lodim} implies that $LO_{A\leftarrow D}\rbag  O_D$ is a dimonoid. Taking into account that $O_D$ and $LO_{A\leftarrow D}$ are rectangular semigroups, we conclude that the dimonoid $LO_{A\leftarrow D}\rbag  O_D$ is rectangular as well. Since in the case $|A|>1$ the semigroup $LO_{A\leftarrow D}$ contains at least two left zeros, and hence it is not commutative, we conclude that the dimonoid $LO_{A\leftarrow D}\rbag  O_D$ is not commutative as well.   Taking into account that in the case $|A|>1$ the semigroup $LO_{A\leftarrow D}$ is not commutative while the semigroup $O_D$ is commutative, we conclude that the semigroups   $LO_{A\leftarrow D}$ and $O_D$ can not be dual, and thus by Proposition~\ref{da} the dimonoid $LO_{A\leftarrow D}\rbag  O_D$ is not abelian.
\end{proof}

\begin{proposition} Let $A$ be a  subset of a set $D$ and $0\in A$. Then  $\Aut(LO_{A\leftarrow D}\rbag  O_D)  = S_{A\setminus\{0\}}\times S_{D\setminus A}$. If $|D|> 2$, then halo of $LO_{A\leftarrow D}\rbag  O_D$ is an emptyset. 
\end{proposition}
\begin{proof} Taking into account that any bijection that preserves $0$ is an automorphism of the null semigroup $O_D$ and according to proof of Proposition~\ref{rcommLOAD} any automorphism of the semigroup $LO_{A\leftarrow D}$ preserves $0$, we conclude that 
$\Aut(LO_{A\leftarrow D}\rbag  O_D) = \Aut(LO_{A\leftarrow D}) = S_{A\setminus\{0\}}\times S_{D\setminus A}$.
Since in the case $|D|> 2$ the commutative semigroup $O_D$  contains no an identity, $\Halo(LO_{A\leftarrow D}\rbag  O_D) = \emptyset$. 
\end{proof}

\begin{remark} By applying the properties of unary duality operations for semigroups and dimonoids discussed in Section~\ref{sec:dual}, we establish that analogous propositions and theorems hold for the semigroups and dimonoids dual to those considered in Sections~\ref{sec:abelian}, \ref{sec:nonabelian}, and \ref{sec:rectangular}.

\end{remark}

\end{document}